\documentstyle{amsart}

\makeatletter
 \def\LaTeX{\leavevmode L\raise.42ex
   \hbox{\kern-.3em\size{\sf@size}{0pt}\selectfont A}\kern-.15em\TeX}
\makeatother

\newcommand{\BibTeX}{{\rm B\kern-.05em{\sc i\kern-.025emb}\kern-.08em\TeX}}

\newtheorem{thm}{Theorem}[section]
\newtheorem{lem}[thm]{Lemma}

\theoremstyle{definition}

\numberwithin{equation}{section}

\begin{document}

\title{Irregular Sampling and the Radon Transform}

\author{Eric Grinberg}
\address{Department of Mathematics, Temple University,
Philadelphia, PA 19122}
\email{grinberg@@math.temple.edu}

\author{Isaac Pesenson}
\address{Department of Mathematics, Temple University,
Philadelphia, PA 19122}
\email{pesenson@@math.temple.edu}

\subjclass{00000; Secondary 00000}

\begin{abstract} In image reconstruction there are techniques that use
analytical formulae for the Radon transform to recover an image from a
continuum of data. In practice, however, one has only discrete data
available. Thus one often resorts to sampling and interpolation methods.
This article presents an approach to the inversion of the Radon
transform that uses a discrete set of samples which need not be
completely regular.

\end{abstract}

\maketitle

\section{Introduction}

The Radon transform of a function $\varphi (x)$ in the plane ${\mathbb R}^d$ is
defined by
$$
R \varphi ( \theta , s) = \int_{x \cdot \theta = s} \, \varphi(x) \, dx,
$$
whenever the integral makes sense. Here $\theta$ is a unit direction
vector and $s$ is a scalar translation parameter. A principal problem in
image reconstruction is the recovery of the values of $\varphi (x)$ from the
data $\{ R\varphi ( \theta , s )\}$ for all $\theta$ and all $s$. The algorithm
that is commonly called {\it Fourier reconstruction} \cite{Nat} is a
discretization of the Fourier-slice or projection-slice formula:

$$
\widehat{\varphi}(\theta
\tau)=(2\pi)^{(1-d)/2}\widehat{(R\varphi)} (\theta,\tau),
$$
where on the right one has the Fourier transform in the second variable
of the Radon transform $R\varphi$.

We consider a finite set of directions $\theta_{j}\in S^{d-1}$, $j=1, 2,
....p$. Using a set of equally spaced samples of the Radon transform $
(R\varphi)(\theta_{j}, s_{\gamma}), \gamma =1, 2,...,q,$ we can
reconstruct $ (R\varphi)(\theta_{j}, s)$ as a function of one variable
$s$. The common way of reconstruction is by applying the
Shannon-Whittaker formula. Taking the Fourier transform in the single
variable $s$ we obtain functions $\widehat{(R\varphi)}(\theta_{j},\tau)$.

Given functions $\widehat{(R\varphi)}(\theta_{j},\tau)$ along all rays
$\theta_{j}\in S^{d-1}$, $j=1, 2, ....p$, we can estimate their values
on a certain polar grid and then reconstruct the function
$\hat{\varphi}(\theta \tau)$. The usual way of reconstruction is again
through the Shannon-Whittaker sampling theorem.

It is well known \cite{Ber}, \cite{Nat} that one of the main problems
with the Fourier reconstruction algorithm is that the Shannon-Whittaker
sampling theorem can be used only in the case of lattice points (regular
sampling). But in many situations there is no way to construct such a
uniform cartesian grid using the naturally available polar grid. A
number of different ways to avoid this obstacle can be found in the book
of Natterer \cite{Nat}. But, in any case, this difficulty reduces the
accuracy of the Fourier reconstruction algorithm. Our idea is to use a
sampling theorem which does not require the uniformity property of the
sample points. The considerations in the present paper are purely
qualitative. The paper \cite{CHM} gives, in the case $d=2$, a different
approach to the reconstruction of the image using an irregular set of
samples.

\section{An irregular sampling theorem}

\noindent
In what follows we use the notations below.

\medskip

\noindent
\newline $B_{\sigma }(\mathbb R^{d}) \, :$ \hfil

\bigskip

\noindent This denotes the set of band limited functions
$B_{\sigma}({\mathbb R}^{d})$ is the set of all $f\in L_{2}({\mathbb R}^{d})$
such that the Fourier transform

$$
\hat{f}(\xi)=(2\pi)^{-d/2}\int_{{\mathbb R}^{d}}f(x)e^{-ix\xi}dx
$$
has support in the ball $B(0,\sigma)$ of radius $\sigma$ centered at
$0$.

\bigskip

\noindent
\newline $\|f\| \,:$ \hfil

\bigskip

\noindent
The symbol $\|f\|$ denotes the $L_{2}(R^{d})$- norm of $f$.

\bigskip\noindent
$X(\lambda ) \, :$ \hfil
\bigskip

\noindent Let $\lambda$ be a positive number. Let $X(\lambda)$ denote a
countable set of points $\{x_{\gamma}\}$ in ${\mathbb R}^d$ with the
following property.
\medskip

There exist a system of open sets $Q(x_{\gamma}, \lambda) \subset
{\mathbb R}^d$ such that

\begin{itemize}
\item
Each $\{x_{\gamma}\}$ contains exactly one point among the
collection $\{x_{\gamma}\}$.
\item
The closure $\overline{Q(x_\gamma , \lambda)}$ is diffeomorphic, as a
manifold with boundary, to a closed ball.
\item
Each $Q(x_{\gamma}, \lambda )$ is of diameter $\leq\lambda$.
\item
The sets $\{ Q(x_{\gamma} , \lambda )\}$ are pairwise disjoint.
\item
 The closures $\overline{Q(x_{\gamma} , \lambda )}$ cover ${\mathbb R}^d$.
\end{itemize}

\medskip

 We will often call $X( \lambda )$ the {\it knot set}.

\bigskip\bigskip

The classical Shannon-Whittaker sampling theorem says that if $f\in
L^{2}({\mathbb R})$ and its Fourier transform $\hat{f}$ has support in
$[-\omega, \omega]$, then $f$ is completely determined by its values at
points $n\Omega $, where $\Omega =\pi /\omega $ and, in the
$L^{2}$-sense,

$$ f(t)= \sum f(n\Omega )\frac{\sin(\pi (t-n\Omega ))}{\pi (t-n\Omega )}.$$

The functions $f\in L^{2}({\mathbb R})$ with the property supp$\hat{f}
\subset [-\omega, \omega ]$ form the Paley-Wiener class $PW_{\omega }$.
The Paley-Wiener theorem states that $f$ is in $PW_{\omega}$ if and only
if $f$ is an entire function of exponential type $\omega $.

Entire functions of finite exponential type are also uniquely determined
by and can be recovered from their values on specific irregular sets of
points ${x_{n}}$. As was shown by Paley and Wiener it is enough to
assume that the functions $\exp ix_{n}t, n\in {\mathbb Z}$ form a Riesz
basis for $L^{2}([-\pi,\pi ])$.

One can consider even more general assumptions about the sequence
$\{ x_{n}\}$. New and old results in the case when the functions $\exp ix_{n}t
$ form different kinds of frames in $L^{2}( [-\omega,\omega])$ were
summarized in \cite{Ben} and \cite{FG}.

Our goal is to show that every band limited function can be
reconstructed from an appropriate irregular set of points using
translations of the fundamental solution of any operator of the form
$\Delta+\varepsilon,\varepsilon\geq 0$, where $\Delta$ is the Laplacian
in Euclidean space. A similar result for the operator $\Delta$ was
considered in \cite{Pes}. We consider the operator
$D=D_{\varepsilon}=\Delta+\varepsilon$, $\varepsilon\geq 0$. The
fundamental solution $E^{k}=E_{\varepsilon}^{k}$ of the operator $D^{k}$
is the inverse Fourier transform of the function
$(|\xi|^{2}+\varepsilon^{2})^{-k}$. In the case when $k>d/2$ and
$\varepsilon>0$ this function is smooth in $L_{2}({\mathbb R}^{d})$ and has fast
decay at infinity (see below). The last property illustrates an
important difference between the cases $\varepsilon=0$ and $\varepsilon
>0$.

If $k>d/2$ and $ \varepsilon>0$ then $(|\xi|^{2}+\varepsilon)^{-k}$ is an
integrable function and, because it is radial, its Fourier transform can
be expressed in terms of the Bessel functions $J_{d/2-1}$:

$$
E^{k}_{\varepsilon}(x)=(2\pi)^{-d}\int_{{\mathbb R}^{d}}\frac{e^{i\xi x}}
 {(|\xi|^{2}+\varepsilon)^{k}}
 d\xi=\frac{|x|^{1-d/2}}{(2\pi)^{d/2}}\int^{\infty}_{0}
\frac{t^{d/2}}{(t^{2}+\varepsilon)
 ^{k}}J_{d/2-1}(t|x|)dt=
$$

$$
\frac{2^{1-k}}{(2\pi)^{d/2}\Gamma(k)}(\varepsilon^{-1/2}
 |x|)^{k-d/2}K_{k-d/2}(\varepsilon^{1/2}|x|).
$$

The last function is a locally integrable since $K_{\nu}(t)$ grows as
$t^{-|\nu|}\ln t$ for $t\rightarrow 0$. It is also a rapidly decreasing
$C^{\infty}$ function outside the origin. We show in the Lemmas 2.2 and
2.3 that for every $k\geq d/2$ there are infinite linear combinations
$L_{x_{\nu}}^{k}\in L_{2}(R^{d})$ of translates of the fundamental
solution $E^{k}$ of the operator $D^{k}$ for which
$L_{x_{\nu}}^{k}(x_{\gamma})= \delta_{\gamma,\nu}, \, x_{\gamma}\in
X(\lambda).$

In general one does not know an explicit formula for $L_{x_{\nu}}^{k}$ so
one has to find approximations to $L_{x_{\nu}}^{k}$ using finite sets of
knots. It is possible to do so because the fundamental solutions
$E^{k}_{\varepsilon}, k>d/2, \varepsilon>0$ have fast decay at infinity.
The value of $L_{x_{\nu}}^{k}$ at a point depends essentially on a
finite number of points from the knot set $X(\lambda)$. In what follows
we will use the notation $L_{\nu}^{k}$ for $L_{x_{\nu}}^{k}$.

\noindent
We prove the following.

\begin{thm} There exists a constant $c=c(d, \epsilon )$ that depends
only on the dimension $d$ and the parameter $\varepsilon$ such that for
any $\sigma >0$ every knot set $ X(\lambda)$ with
$\lambda<((\sigma+\varepsilon) c(d, \varepsilon))^{-1} $ and every
integer $r\geq[d/2]+1$,

$$
f= \lim_{l\rightarrow\infty}\sum_{x_{\nu}\in X(\lambda)}f(x_{\nu})
L^{2^{l}r}_{\nu}, \, l\in {\mathbb N},
$$
for all $f\in B_{\sigma}(R^{d})$.

\noindent
Moreover, an error estimate for this approximation is

$$\|f-\sum_{\nu}f(x_{\nu})L^{2^{l}r}_{\nu}
\|\leq
 2(c(d, \varepsilon)\lambda (\sigma +\varepsilon))^{2^{l+1}r}\|f\|.$$

\end{thm}

\medskip\noindent
The proof of the Theorem will follow from some preliminary results
below.

\bigskip\bigskip

Given a sample set $X(\lambda)$ and a sequence $\{s_{\gamma}\}\in l_{2}$
we will be interested in finding a function $s_{k}\in H^{2k},$ for $k$
large enough, such that
\bigskip

a) $ s_{k}(x_{\gamma})=s_{\gamma},$ for $ x_{\gamma}\in \ X(\lambda).$

b) The function $s_{k}$ minimizes the functional $u\rightarrow \|D^{k}u\|$.

\bigskip

In what follows we  will use the fact that, in the case $\varepsilon>0$,
the functional $u\rightarrow \|D^{k}u\|$ is equivalent to the Sobolev
norm.

For the given sequence $\{s_{\gamma} \}\in l_{2}$ consider a function $f$
from $H^{2k}$ such that $f(x_{\gamma})=s_{\gamma}.$ Let $Pf$ denote the
orthogonal projection of this function $f$ (in the Hilbert space
$H^{2k}$ with the natural inner product) on the subspace

$$U^{2k}(X(\lambda))=\left \{f\in H^{2k}|f(x_{\gamma})=0\right \}
$$
with $H^{2k}$-norm. Then the function $g=f-Pf$ will be the unique
solution of the above minimization problem for the functional
$u\rightarrow \|D^{k}u\|, k>d/2$.

Given a function $f\in H^{k},$ where $ k>d/2 $, the function $s_{k}(f)$
will denote the solution to the above optimization problem with
$s_{\gamma}=f(x_{\gamma})$.

\bigskip We will denote by $S^{2k}(X(\lambda))$ the set of all $L_{2}$-
solutions of the equation

$$
D ^{2k}u=\sum_{x_{\gamma}\in X(\lambda)}\alpha _{\gamma }\delta (x_{\gamma
}),
$$
where $\delta (x) $ is the Dirac measure and $\{ \alpha _{\gamma} \} \in
l_{2}$. Our next goal is to show that every $s_{k}(f)$ belongs to $
S^{2k}(X(\lambda))$.

 Indeed, suppose that $s_{k}\in H^{2k}$ is a solution to the
minimization problem and $h\in U^{2k}(X(\lambda)).$ Then

$$\|D ^{k}(s_{k}+\lambda h)\|^{2}=\|D^{k}s_{k}\|^{2}_{2}+2Re
{\lambda\int
 D^{k}s_{k} D^{k}h}dx+|\lambda |^{2}\|D ^{k}h\|^{2}_{2}.$$

The function $s_{k}$ can be a minimizer only if for any $h\in
U^{2k}(X(\lambda))$

$$
\int D^{k}s_{k}D^{k}hdx=0.
$$

So, the function $g= D^{k}s_{k}\in L_{2} $ is orthogonal to
$D^{k}U^{2k}(X(\lambda)) $. Let $\varphi_{\gamma}\in C^{\infty}_{0}$
have disjoint supports and $\varphi_{\gamma}(x_{\gamma})=1$ and $h\in
C_{0}^{\infty}$. Then the function $h-\sum
h(x_{\gamma})\varphi_{\gamma}$ belongs to the space $U^{2k}
(X(\lambda))\cap C_{0}^{\infty}.$ Thus,

$$
0=\int g \overline{D ^{k}(h-\sum
h_{\gamma}\varphi_{\gamma})}dx=\int g \overline{D^{k} h}dx-
\sum \overline{h(x_{\gamma})}\int g \overline{D^{k}\varphi_{\gamma}}dx.
$$

In other words

$$
D ^{k}g =\sum_{x_{\gamma}\in X(\lambda)}\alpha _{\gamma }\delta
(x_{\gamma }),
$$
or
$$
D^{2k}s_{k}=\sum_{x_{\gamma}\in X(\lambda)}\alpha _{\gamma }\delta
(x_{\gamma }),
$$
where $\delta (x) $ is the Dirac measure.

Moreover, for any integer $r>0$

$$
\sum _{\gamma =1}^{r}|\alpha _{\gamma }|^{2}=<\sum_{1}^{\infty}\alpha
_{\gamma
 } \delta (x_{\gamma}), \sum _{1}^{r}\alpha _{\gamma} \varphi _{\gamma}>
\leq C \|\sum_{1}^{\infty} \alpha _{\gamma} \delta
(x_{\gamma})\|_{H^{-2k}} (\sum_{1}^{r} |\alpha
_{\gamma}|^{2})^{1/2},
$$
where $C$ is independent of $r$. This shows that the sequence $\{\alpha
_{\gamma }\}$ belongs to $l_{2}$.

 Now suppose that $f\in H^{\infty}$ and

$$D ^{2k}f=\sum_{x_{\gamma}\in X(\lambda)}\alpha _{\gamma }\delta
(x_{\gamma}),
$$

 where $\{ \alpha _{\gamma} \} \in l_{2}$.

It was shown in \cite{Pes} the norm of the Sobolev space $H^{k}$ is
equivalent to the norm

$$\|D^{k/2}f\|+(\sum|f(x_{\gamma})|^{2})^{1/2}.$$

So for any $\mu >0$ we have

$$|<D ^{2k}f , g>|= | <\sum \alpha _{\gamma}\delta(x_{\gamma}), g>| \leq
$$ $$ \left(\sum |\alpha _{\gamma}|^{2}\right)^{1/2}\left(\sum
|g(x_{\gamma})^{2}\right)^{1/2} \leq C\left(\sum |\alpha
_{\gamma}|^{2}\right)^{1/2}\|g\|_{H^{d/2+\mu }}.$$ This shows that the
distribution $\sum _{1}^{\infty}\alpha _{\gamma} \delta
(x_{\gamma})=D^{2k} f$ belongs to $H^{-d/2-\mu}$. Since the operator $D
^{2k}$ is $C^{\infty }-$uniformly elliptic of order $ 4k$ we can use a
corresponding regularity result which gives that $f$ belongs to
$H^{-d/2-\mu+4k}$, which is included in $H^{2k}$ for all $k>d$. The
assertion that the orthogonal complement of $D ^{k} U^{2k}(X(\lambda)) $
is a subset of $ S^{2k}(X(\lambda))$ is proved.

 Conversely, if $f , h $ belong to $ S^{2k}(X(\lambda))$ and
$U^{2k}(X(\lambda))
\cap
C^{\infty}$ respectively, then, since
 $f\in H^{2k}$ and $h\in H^{2k}$ and the pairing $< .  , .  >$ is an
extension of the scalar product in $L^{2}$,

$$
\int\overline{D ^{k}h}dx =<D^{k}f, \overline{h}>=
\sum \alpha_{\gamma}\overline{h(x_{\gamma})}=0.$$

Thus we proved the following.

\begin{lem}
 A function $f\in L_{2}$ belongs to  $ S^{2k}(X(\lambda))$, i.e.
satisfies the equation

$$
D ^{2k}f=\sum_{x_{\gamma}\in X(\lambda)}\alpha _{\gamma }
\delta (x_{\gamma }),
$$

\noindent
where $\{ \alpha _{\gamma} \} \in l_{2}$ if and only if $f$ is a
solution to the minimization problem for the functional
$u\rightarrow \|D^{k}u\|$.

In particular, every solution to the minimization problem is a linear
combination of the translates of the fundamental solution
$E^{k}_{\varepsilon}$ of the operator $D_{\varepsilon}^{k}
=(\Delta+\varepsilon)^{k}, \varepsilon>0$. \end{lem}
 In particular, for any $x_{\gamma}\in X(\lambda)$ there exists a unique
$L^{2k}_{\gamma}(X(\lambda))\in S^{2k}(X(\lambda) ) $ that takes
the value $1$ at the point $x_{\gamma}$ and $0$ at all other points in
$X(\lambda)$. These functions form a Riesz basis in
$S^{2k}(X(\lambda)).$

\medskip

Recall that the last assertion means that for any
$g\in S^{2k}(X(\lambda))$ in $L_{2}$ we have
$$
g=\sum_{\gamma}g(x_{\gamma})L^{2k}_{\gamma}
$$
and there are constants $C_{1},C_{2}>0$ such that
$$
\|g\|_{2}\leq C_{1}\left( \sum
|g(x_{\gamma})|^{2}\right)^{1/2}\leq C_{2}\|g\|_{2}, \qquad (g\in
S^{2k}(X(\lambda))).
$$

\medskip

 This statement is a consequence of the next result.

\begin{lem} Every function from $S^{2k}(X(\lambda)), \, k=2^{l}d$, is
uniquely determined by its values at points $x_{\gamma}\in X(\lambda)$.
Moreover, for any $f\in S^{2k}(X(\lambda) ) $ the norm
$(\sum |f(x_{\gamma})|^{2})^{1/2}$ is equivalent to the $L_{2}$-norm
and to the Sobolev norm. \end{lem}

\medskip

 \begin{pf}
Since $S^{2k}(X(\lambda))$ is closed in the $L_{2}$-norm and
$S^{2k}(X(\lambda))\subset H^{2k}$ the $L_{2}$-norm and $H^{2k}$ norm
are equivalent on $S^{2k}(X(\lambda))$. Moreover, one can show that on
the space $S^{2k}(X(\lambda)), k=2^{l}d,$ the norm $H^{2k}$ is
equivalent to the norm
 $(\sum |f(x_{\gamma})|^{2})^{1/2}, f\in S^{2k}(X(\lambda))$.

Indeed, if the functions $\varphi _{\gamma }\in C^{\infty}$ have
disjoint supports in $B(x_{\gamma},\lambda /4)$ and $\varphi
_{\gamma}(x_{\mu})=\delta _{\gamma \mu }, |\varphi _{\gamma}|\leq 1$,
then the function $F=\sum_{\gamma \in N} f(x_{\gamma })
\varphi_{\gamma}$ is in $H^{2k}$ and $f(x_{\gamma})=F(x_{\gamma}), k >
d/2$. Because of the minimization property, we have

$$
\|D ^{k}f\|\leq \|D^{k}F \|\leq C\left(\sum_{\gamma}
|f(x_{\gamma})|^{2}\right)^{1/2}.
$$
Since for $k=2^{l}d$ the $H^{2k}$ norm on $S^{2k}$ is equivalent to
the norm
$\|D ^{k}f\|$,  this
 implies its equivalence to the norm $\left( \sum_{\gamma}
|f(x_{\gamma})|^{2}\right)^{1/2}.$

\end{pf}

Now we can prove the following approximation property.

\begin{thm} For any integer $r\geq[d/2]+1$ and any $f\in
H^{2^{l+1}r}(R^{d})$,

$$
f(x)=\lim_{l\rightarrow\infty}s_{2^{l}r}(f)
    =\lim_{l\rightarrow\infty} \sum_{x_{\nu}\in X(\lambda)}
 f(x_{\nu}) L^{2^{l}r}_{\nu}(x).
$$
Moreover, there exists a constant $c(d,\epsilon )$ that depends  only on the dimension
$d$ and the parameter $\varepsilon$ such that the following error estimate
is valid:
$$
\|f-\sum_{\nu}f(x_{\nu})L^{2^{l}r}_{\nu})\|\leq
 2(c(d, \varepsilon)\lambda)^{2^{l+1}r}\|D^{2^{l}r}f\|, \, l=0, 1, \ldots
$$

\end{thm}

\begin{pf}
If
$$
f\in H^{2k}, k=2^{l}d, \qquad l=0,1,...
$$
and
$$
s_{k}(f)= \sum_{\nu}f(x_{\nu})L^{k}_{\nu}
$$
then
$$
f-s_{k}(f)\in U^{2k}(X(\lambda))
$$
and as it was shown in \cite{Pes} we have

$$
\|f-s_{k}(f)\|\leq
 (C(d,\varepsilon) \lambda)^{k}\|D^{k/2}(f-s_{k}(f))\|, k=2^{l}d,
(l=0,1, \ldots ) .
$$

Using the minimization property of $s_{k}(f)$ we obtain

$$
\|f-s_{k}(f)\|\leq (c(d, \varepsilon)\lambda )^{k}\|D
^{k/2}f\|,c(d,\varepsilon)=
2C(d,\varepsilon), k=2^{l}d, \qquad
(l=0,1, \ldots ) .
$$
The approximation theorem is proved.
\end{pf}

Our Theorem 2.1 follows from the above approximation theorem and the
Bernstein inequality satisfied by any function from $B_{\sigma}$:

$$\|D^{k/2}f\|\leq(\sigma +\varepsilon)^{k}\|f\|.$$

Theorem 2.1 is proved.

\bigskip

% remarks Although in this paper we consider primarily the case of an
irregular set of knots and the operator $\Delta +\varepsilon$, for
$\varepsilon>0$, the case of equally spaced points for the operator
$\Delta$ is of special interest because in this case explicit formulas
for the Fourier transform of $L_{\bar{n}}^{k}$ are known; here $\bar{n}$
is the integer lattice. Indeed one can verify (see \cite{Mad}) that
in the case of the standard lattice $\bar{n}$ of ${\mathbb R}^{d}$ the function
$\Lambda^{k}_{0}=\hat{L^{k}_{0}}$ is

$$
\Lambda^{k}_{0}(\xi)=(2\pi)^{-d/2}
(|\xi|^{2k}\sum_{\bar j\in {\mathbb Z}^{d}}|\xi-2\pi\bar j|^{-2k})^{-1}
$$
and all other $L^{k}_{\bar{n}}$ are translations of $L^{k}_{0}$.

These functions $L_{\bar{n}}^{k}$ have very fast decay at infinity in
the sense that for every $k$ there are $a=a(k)>0, b=b(k)>0$ such that
$|L_{\bar{n}}^{k}(x)|\leq ae^{-b|\bar{n}-x|}$.

We also want to make the following remark. Our sampling theorem requires
in general some {\it oversampling}. This means that the distance between
sampling points needs to be small enough compared to the size of the
support of the Fourier transform of the given band limited function. We
will show now that if one is going to consider lattice sampling points
then the oversampling is not necessary. For example if the Fourier
transform of a function is in the cube $[-\pi,\pi]^{d}$ then the natural
lattice in ${\mathbb R}^{d}$ can be chosen as the sampling set. Such a
rate of sampling is known to be the best possible and is called the {\it
Nyquist rate}.

Indeed, we can rewrite the formula for the function $\Lambda^{k}_{0}$:

$$
\Lambda^{k}_{0}(\xi)=(|\xi|^{-2k}|\xi-2\pi\bar{j}|)^{2k}
\Lambda^{k}_{2\pi\bar{j}}(\xi).
$$

This shows that $\lim_{k\rightarrow\infty}\Lambda^{k}_{0}(\xi)$ is zero
for every $\xi$ outside the cube $[-\pi,\pi]^{d}.$

Next, the formula

$$
\Lambda^{k}_{0}(\xi)=(2\pi)^{-d/2}(1+|\xi|^{2k}\sum_{\bar{j}\in {\mathbb Z}^{d}_{+}}
|\xi-2\pi\bar{j}|^{-2k})^{-1},
$$
where ${\mathbb Z}^{d}_{+}$ is the set of all non zero $d$-tuples, implies that
the limit $\lim_{k\rightarrow\infty}\Lambda^{k}_{0}(\xi)$ is
$(2\pi)^{-d/2}$
for all $\xi $ in the cube $[-\pi,\pi]^{d}.$

 In other words

$$
\lim_{k\rightarrow\infty}L^{k}_{0}(x)=
\frac{\sin(\pi x_{1})}{\pi x_{1}}\frac{\sin(\pi x_{2})}{\pi x_{2}}...
\frac{\sin(\pi x_{d})}{\pi x_{d}}.
$$

 Together with the classical Shannon-Whittaker sampling theorem,

$$
\hat{\varphi}(t)= \sum_{\bar{n}\in Z^{d}} \hat{\varphi}(\bar{n} )
\frac{\sin(\pi (t_{1}-n_{1} ))}
{\pi (t_{1}-n_{1} )} \ldots
\frac{\sin(\pi (t_{d}-n_{d} ))}{\pi (t_{d}-n_{d} )},
$$
where $\varphi$ has support in the cube $[-\pi,\pi]^{d}$,
 $\bar{n}=(n_{1}, \ldots ,n_{d}), t=(t_{1}, \ldots ,t_{d}),$ this proves
the formula

$$
\varphi(x)=\lim_{k\rightarrow\infty}\sum_{\bar{n}\in Z^{d}}
\hat{\varphi}(\bar{n})\Lambda^{k}_{\bar{n}}(x), \qquad  x\in {\mathbb R}^{d}.
$$

 This interpolation formula seems to be new.

\bigskip\bigskip

\section {Inversion of the Radon transform in ${\mathbb R}^{n}$ using irregular
sampling}

\bigskip

For a given function $\varphi$ on ${\mathbb R}^{d}$ the Radon transform
$R_\theta \varphi$ is defined by

$$
R_{\theta}\varphi(s)=\int_{\theta ^{\perp}}\varphi(s\theta +y)dy,
$$
\noindent
where $\theta $ is a direction vector belonging to the unit sphere
$S^{d-1}$ and $s$ is a real number. In other words the Radon transform
$R_{\theta}\varphi(s)=R\varphi(\theta, s)$ is the integral of $\varphi$
over the hyperplane in ${\mathbb R}^{d}$ defined by ${x: <x,\theta>=s}$. The
backprojection operator is defined by

$$ R^{\ast}g(x)=\int_{S^{d-1}}g(\theta, <x, \theta>)d\theta, $$ where
$x\in R^{d}$, $g$ is defined on the direct product of $S^{d-1}$, and
$\mathbb R$, which can be identified with the set of hyperplanes in
${\mathbb R}^{d-1}$.

Then, if $\varphi \in C^{\infty}_{0}(R^{d})$, the identity

$$
R^{\ast}I^{1-d}R\varphi=\varphi,
$$
holds, where, for $\alpha$ is real and $ I^{\alpha}$ is the Riesz potential
operator, i.e., the Fourier transform of the function $I^{\alpha}\varphi$
is defined as $|\xi|^{-\alpha}\hat{\varphi}(\xi)$. For proofs see \cite{Nat}.
Our goal is to introduce a different reconstruction formula which only
requires a discrete set of values of the Radon transform.

The analogous formula in Fourier analysis is the Poisson summation
formula for the functions $\varphi$ from $L_{2}(R)$ with support in
$[-\pi,\pi]$:

\begin{equation}\varphi(t)=\sum_{n\in
Z}\hat{\varphi}(n)e^{int}. \end{equation}

The meaning of the last formula is that the Fourier coefficients of a
function with compact support are regularly spaced samples of its
Fourier transform. In this paper we will give an analog of the Poisson
summation formula for the Radon transform. More precisely it will be
shown that a compactly supported function can be reconstructed using
even an irregular set of samples of its Radon transform.

Applying the Fourier transform to the formulas from the Theorem 2.1 we
arrive at the following irregular version of the Poisson summation
formula (3.1).

\begin{thm}
If $\varphi\in L_{2}(R^{d})$ has support in the ball $B(\sigma, 0)$ and
$\Lambda^{k}_{\nu}$ is the inverse Fourier transform of the function
$L^{k}_{\nu}$ then

\begin{equation}
\varphi = \lim_{l\rightarrow\infty}\sum_{\xi_{\nu}\in\Xi(\lambda)}
\hat{\varphi}(\xi_{\nu})
\Lambda^{2^{l}r}_{\nu},\end{equation}
assuming $l\in \mathbb N$, $r\geq[d/2]+1$,
$\lambda<(c(d, \varepsilon)(\sigma+\varepsilon))^{-1}$,
and  where the $\Xi(\lambda)$ is an appropriate discrete set in
the space of the dual variable $\xi$.
\hfil \newline
 An error estimate for this approximation is

$$
\|\varphi-\sum_{\xi_{\nu}\in\Xi(\lambda)}\hat{\varphi}(\xi_{\nu})\Lambda^{2^{l
}r}_{\nu}
\| \leq 2(c(d, \varepsilon)\lambda
(\sigma+\varepsilon)^{2^{l+1}r}\|{\varphi}\|.
$$
 \end{thm}

\bigskip\bigskip
 \noindent
 Note that  $\Lambda^{k}_{\nu}=\Lambda^{k}_{\nu, \varepsilon}$ is of the form

\begin{equation}
 \Lambda^{k}_{\nu, \varepsilon}(x)=
(|x|^ {2}+\varepsilon)^{-k}\sum_{\xi_{\mu}\in\Xi(\lambda)}
 a_{\mu}(\nu, k)
  \exp ( -i\xi_{\mu}x)
 \end{equation}
 and the coefficients $ a_{\mu}(\nu, k)$ can be determined using the conditions
 $L^{k}_{x_{\nu}}(x_{\sigma})
 =\delta_{\nu, \sigma} $, where

 $$L^{k}_{\nu, \varepsilon}(x)=\sum_{\mu}a_{\mu}(\nu, k)
 E^{k}_{\varepsilon}(x-x_{\mu}).$$

\medskip\noindent
The following statement is an analog of the last theorem in the case of
the Radon transform.

\begin{lem} Suppose that $\varphi \in C_{0}^{\infty}$ has support in the
ball $B(\sigma, 0)$. If $\Xi(\lambda)$ is a knot set in the space of the
dual variable $\xi$ with $\lambda<(c(d,
\varepsilon)(\sigma+\varepsilon))^{-1},$ where $c(d, \varepsilon)$ is
from the Theorem 2.1, then

 \begin{equation}
\varphi(x)=\lim_{l\rightarrow\infty}(2\pi)^{(1-d)/2}\sum
_{x_{\nu}  \in\Xi(\lambda)}
\widehat {(R\varphi)}(\xi_{\nu})\Lambda^{2^{l}r}_{\nu}(x) \quad
(l\in {\mathbb N}, x\in {\mathbb R}^{d}),\end{equation}
where $\Lambda^{k}_{\nu}$ are from (3.3).
\newline
An error estimate is given by the inequality

$$\|\varphi-(2\pi)^{(1-d)/2}\sum_{\nu}\widehat{(R\varphi)}(\xi_{\nu})
\Lambda^{2^{l}r}_{\nu}\|
 \leq
 2(c(d, \varepsilon)\lambda (\sigma+\varepsilon))^{2^{l+1}r}\|R\varphi\|.$$

 \end{lem}

\begin{pf}
 The Fourier slice theorem states that

$$  \hat{\varphi}(\theta \tau)=(2\pi)^{(1-d)/2}\widehat{(R\varphi)}(\theta,
\tau),$$
where on the right we have the Fourier transform in the second variable.

Let $\varphi \in L_{2}({\mathbb R}^{d})$ be supported in the ball $B(0,\sigma)$.
 Then the function $f$ such that $f=\hat{\varphi}$ belongs to
 $B_{\sigma}({\mathbb R}^{d})$. According to our sampling theorem,

$$
f(\xi)=\lim_{l\rightarrow\infty}\sum_{\xi_{\nu}\in \Xi(\lambda)}
 f(\xi_{\nu})L^{2^{l}r}_{\nu}(\xi),
$$
where the $\Xi(\lambda)$ is an appropriate discrete set in
the space of the dual variable $\xi$.
  The error estimate is given by the inequality

$$\|f-\sum_{\nu}f(\xi_{\nu})
L^{2^{l}r}_{\nu}\|\leq
 2(c(d, \varepsilon)\lambda( \sigma+\varepsilon))^{2^{l+1}r}\|f\|.$$

So we have

$$
f(\xi)=\lim_{l\rightarrow\infty}(2\pi)^{(1-d)/2}\sum _{\nu}
\widehat{(R\varphi)}(\xi_{\nu})L^{2^{l}r}_{\nu}(\xi),
$$
where convergence is understood in the $L_{2}({\mathbb R}^{d})$ sense.

Taking inverse Fourier transform we obtain
\begin{equation}\varphi(x)=\lim_{l\rightarrow\infty}(2\pi)^{(1-d)/2}\sum
_{\nu}
\widehat {(R\varphi)}(\xi_{\nu})\Lambda^{2^{l}r}_{\nu}(x) \quad
( l\in N, x\in {\mathbb R}^{d},\end{equation}
where $\Lambda^{k}_{\nu}$
is the inverse Fourier transform of $L^{k}_{\nu}$.
The error estimate is given by the inequality

$$\|\varphi-(2\pi)^{(1-d)/2}\sum_{\nu}\widehat{(R\varphi)}(\xi_{\nu})
\Lambda^{2^{l}r}_{\nu}\|
 \leq
 2(c(d, \varepsilon)\lambda (\sigma+\varepsilon))^{2^{l+1}r}\|R\varphi\|.$$
 \end{pf}

Although the formula (3.4) is a natural analog of the formula (3.2) it
involves not only the Radon transform but also the Fourier transform in
the second variable. We are going to show how one can approximate the
values $\widehat {(R\varphi)}(\xi_{\nu})$ using just samples of the Radon
transform. One way to do this is by using equally spaced samples and the
Shannon-Whittaker formula. This method has the advantage that the Fast
Fourier Transform can be used \cite{Nat}. But it is also available in the
context of our sampling Theorem 2.1 using an irregular set of samples.

It is convenient to assume now that every point in $\Xi(\lambda)$ has
polar coordinates $(\theta_{\nu}, \tau_{\mu})$ where $\theta_{\nu}$
belongs to the unit sphere $S^{d-1}$ and $\tau_{\mu}$ is the distance
from $0$.

In the one dimensional case, $d=1$, we will use the notation $U^{k}_{\nu}$
for the functions $L^{k}_{\nu}$ which were constructed in the second
section. Note that $U^{k}_{\nu}$ is a piecewise polynomial spline of
order $2k$ with knot sequence $x_{\nu}$.

\begin{lem} If $Y(\mu_{j})$ is a sequence of knots with
$\mu_{j}\rightarrow 0$ and $U^{2^{m}}_{\gamma, j}$ is the corresponding
set of one-dimensional Lagrangian splines then

 \begin{equation}\widehat{(R\varphi)}(\theta,\tau)=\lim_{j\rightarrow\infty}
 \sum_{s_{\gamma}\in Y(\mu_{j})}(R\varphi)(\theta,s_{\gamma})
 V^{2^{m}}_{\gamma, j}(\tau)\end{equation}
 where $V^{2^{m}}_{\gamma, j}$ is the Fourier transform of
 $U^{2^{m}}_{\gamma, j}.$

 \end{lem}
 \begin{pf}
 Let $\lambda_{i}\rightarrow 0$. Now  $X_{i}=X(\lambda_{i})$ is the
 corresponding knot set and $\Lambda^{2^{l}r}_{\nu}(X(\lambda_{i}))$
 is the set of Lagrange functions that correspond to $X_{i}.$

 \medskip\noindent
 On the Fourier transform side our approximation theorem gives

 \begin{equation}\hat{g}=\lim_{i\rightarrow\infty} \sum_{x_{\nu}\in
X(\lambda_{i})}
 g(x_{\nu})\Lambda^{2^{l}r}_{\nu}(X(\lambda_{i})), g\in C_{0}^{\infty},
 \end{equation}
 with an error estimate
$$
\|\hat{g}-\sum_{\nu}g(x_{\nu})\Lambda^{2^{l}r}_{\nu}(X(\lambda_{i}))\|\leq
 2(c\lambda_{i})^{2^{l+1}r}\|g\|_{H^{2^{l+1}r}}, \quad (l=0, 1, \ldots ).
$$

Note that the sum (3.6) is finite because $g$ has compact support. To
describe the corresponding approximation to $R\varphi(\theta,s)$ as the
function of the variable $s$ we introduce the sequence of $\mu_{j}
\rightarrow 0$ and corresponding knot sequences $Y_{j}=Y(\mu_{j}).$ Then,
since we are considering the one-dimensional case we have

 $$(R\varphi)(\theta,s)=\lim_{j\rightarrow\infty}
 \sum_{s_{\gamma}\in Y(\mu_{j})}(R\varphi)(\theta,s_{\gamma})
 U^{2^{m}}_{\gamma, j}(s), \, m>1,
$$
and taking the Fourier transform in $s$ we obtain

 $$
\widehat{(R\varphi)}(\theta,\tau)=\lim_{j\rightarrow\infty}
\sum_{s_{\gamma}\in Y(\mu_{j})}(R\varphi)(\theta,s_{\gamma})
V^{2^{m}}_{\gamma, j}(\tau)
$$
where $V^{2^{m}}_{\gamma, j}$ is the Fourier transform of
$U^{2^{m}}_{\gamma, j}.$

 Note that if  $\varphi$ has compact support then its Radon transform
$R_{\theta}\varphi(s)$ is also of compact support in the variable $s$
and, because of this, the last two sums are finite.
 \end{pf}

 Keeping the same notations we summarize the last two lemmas in the
following theorem.

 \begin{thm}
 Suppose that $\varphi \in C_{0}^{\infty}$ has support in the ball
$B(\sigma, 0)
$.
 If $\Xi(\lambda)$ is a knot set in the space of the dual variable $\xi$ with
 $\lambda<(c(d, \varepsilon)(\sigma+\varepsilon)^{-1}$ where $c(d,
\varepsilon)$
 is from the Theorem 2.1 then

 \begin{equation}\varphi(x)=\lim_{l\rightarrow\infty}(2\pi)^{(1-d)/2}\sum
_{\nu,
\mu}
\widehat {(R\varphi)}(\theta_{\nu},\tau_{\mu})\Lambda^{2^{l}r}_{\nu,\mu}(x),
l\in N \quad (x\in {\mathbb R}^{d}),\end{equation}
where $\Lambda^{k}_{\nu,\mu}$
is the inverse Fourier transform of $L^{k}_{\nu,\mu}$ and of the form (3.3).

An error estimate is given by the inequality

$$\|\varphi-(2\pi)^{(1-d)/2}\sum_{\nu,
\mu}\widehat{(R\varphi)}(\theta_{\nu},\tau_{\mu})
\Lambda^{2^{l}r}_{\nu,\mu}\|
 \leq
 2(c(d, \varepsilon)\lambda (\sigma+\varepsilon))^{2^{l+1}r}\|R\varphi\|.
$$

The approximate values of $\widehat{R\varphi}(\theta_{\nu},\tau_{\mu})$
can be determined by the formula

 $$
\widehat{(R\varphi)}(\theta_{\nu},\tau_{\mu})=\lim_{j\rightarrow\infty}
 \sum_{s_{\gamma}\in Y(\mu_{j})}(R\varphi)(\theta_{\nu},s_{\gamma})
 V^{2^{m}}_{\gamma, j}(\tau_{\mu}).
$$

\end{thm}

\bigskip\bigskip

 \section{A Computational algorithm}

\bigskip\bigskip

Alas, the amount of work which is needed for numerical implementations
of the above algorithm is too big. In what follows we improve the
standard Fourier Algorithm by using our irregular sampling theorem in
conjunction with the Fast Fourier Transform. It will make our
modification as efficient as the original Fourier Algorithm is. We
restrict ourselves to the case $d=2$ and $\varepsilon=1$. In this case
the constant $c(d,\varepsilon)$ is not greater than $3$.

Let $\varphi $ have compact support and let $g=R\varphi$ be sampled at
$$
(\theta_{j},s_{l}), \quad j=1, \ldots ,p, \qquad l=-q, \ldots ,q,
$$
where
$$
\theta_{j}=(\cos\alpha_{j}, \sin\alpha_{j}), \quad  \alpha_{j}=(j-1)\pi/p, s_{l}=hl, h=1/q.
$$

It is easy to see that the optimal relation between $p$ and $q$
is given by the approximate formula $p \approx \pi q$ .

\bigskip\bigskip

\noindent
{\bf STEP 1.} For $j=1,...,p,$ compute approximations $\hat{g}_{jr}$ to $\hat{g}
 (\theta_{j}, r\pi)$ by

$$
\hat{g}_{jr}=(2\pi)^{-1/2}h\sum_{l=-q}^{q-1}e^{ -i\pi lr/q }g(\theta_{j},
 s_{l}), \, r=-q,...,q-1.
$$

 This step provides an approximation to $\hat{\varphi}$ on the polar grid

 $$
G_{p, q}={\pi r \theta_{j} :r= -q,...,q-1, j= 1,...,p}.
$$

 We have $\hat{\varphi}(r\pi \theta_{j})=(2\pi)^{-1/2}\hat{g}_{rj}$.

 Because we perform $p$ discrete Fourier transforms of length $2q$ and
$p=\pi q$ this step requires $0(q^{2}\ln q)$ operations.

\bigskip

\noindent
{\bf  STEP 2.} For each $k\in {\mathbb Z}^{2}, |k|<q$ use the formula

$$
\hat{\varphi}_{k}= \sum \hat{g}_{jl} L^{m}_{jl}
$$
where the summation is taken over some points from the polar grid
$G_{pq}$ which surround the point $k\in {\mathbb Z}^{2}$. Moreover due to the fact
that the function $ L^{m}_{jl}$ is localized essentially around point
$\pi l \theta_{j}$ it is enough to keep a constant number of terms in
this summation. This observation is very important since it implies that
the second step requires essentially $0(q^{2})$ operations.

\bigskip

\noindent
{\bf  STEP 3.}
Compute an approximation $\varphi_{N}$ to $\varphi(hN), N\in {\mathbb Z}^{2}$ by

$$
\varphi_{N}=(1/2\pi)^{d/2}\sum_{|k|<q}e^ {i\pi Nk/q}\varphi_{k}, |N|<q.
$$

\bigskip

To perform this step one needs $0(q^{2}\ln q)$ steps. Thus the amount of
work for our modified algorithm is the same as for the standard Fourier
algorithm.

The high stability in the step 2 is a consequence of the estimate from
Theorem 2.4.

 \section{ ACKNOWLEDGMENTS}

 We thank Professor G. Herman for sending us the reprint of the paper \cite{CHM}.

The second author thanks Professors W. Madych and B. Rubin for
stimulating and useful discussions.

\makeatletter \renewcommand{\@biblabel}[1]{\hfill#1.}\makeatother

\end{document}